\documentclass[10pt,a4paper]{article}
\usepackage{amssymb,amsmath}
\usepackage[T1]{fontenc}
\usepackage[english]{babel}
\usepackage{alltt}

\newtheorem{hypo}{Hypothesis}

\newtheorem{prop}[hypo]{Proposition}

\newtheorem{thm}[hypo]{Theorem}

\begin{document}

        \title{A limit theorem for a random walk in a
        stationary scenery coming from a hyperbolic dynamical system}
        \author{Françoise Pène\\
         Universit\'e de Bretagne Occidentale\\
         UMR CNRS 6205\\
         D\'epartement de Math\'ematiques,
         UFR Sciences et Techniques\\
         6, avenue Victor Le Gorgeu, 29238 BREST Cedex 3, France\\
         francoise.pene@univ-brest.fr
}
    \maketitle

\noindent{\bf Abstract.} {\it 
In this paper, we extend a result of Kesten and Spitzer
\cite{KestenSpitzer}. 
Let us consider an invertible probability dynamical system
$(M,{\cal F},\nu,T)$ 
and $f:M\rightarrow{\mathbb R}$
some function with null expectation. 
We define the stationary sequence $(\xi_k:=f\circ T^k)_{k\in\mathbb Z}$.
Let $(S_n)_{n\ge 0}$ be a simple
symmetric random walk on $\mathbb Z$ independent of $(\xi_k)_{k\in\mathbb Z}$. 
We are interested in the study of the 
sequence of random variables of the form
$\left(\sum_{k=1}^n\xi_{S_k}\right)_{n\ge 1}$.
We give examples of partially hyperbolic dynamical systems $(M,{\cal F},\nu,T)$ 
and of functions $f$ such that $\left({1\over n^{3\over 4}}
\sum_{k=1}^n\xi_{S_k}\right)_{n\ge 1}$
converges in distribution.
}

\section{Introduction}
In \cite{KestenSpitzer}, Kesten and Spitzer prove that, if 
$(\xi_n)_{n\in\mathbb Z}$ is a sequence of independent
identically distributed satisfying a central limit theorem
and if $(S_n)_{n\ge 0}$ is the simple symmetric
random walk on $\mathbb Z$ independent of $(\xi_k)_{k\in\mathbb Z}$, 
then $\left({1\over n^{3\over 4}}
\sum_{i=1}^n\xi_{S_k}\right)_{n\ge 1}$ converges in distribution.
In this paper, our goal is to establish such a result
when $(\xi_k)_{k\in\mathbb Z}$ is a stationary sequence of random variables
given by a dynamical system with some hyperbolic properties.
More precisely, we study the cases when $(\xi_k=f\circ T^k)_{k\in\mathbb Z}$,
with $f$ a $\nu$-centered H\"older continuous function 
and when $(M,{\cal F},\nu,T) $ is one of the following dynamical
systems~:
\begin{itemize}
\item the transformation $T$ is
an ergodic algebraic automorphism of the torus $M={\mathbb T}^{d_0}$
endowed with its normalised Haar measure $\nu$ (for some $d_0\ge 2$);
\item  the transformation $T$ is a diagonal transformation 
on a compact quotient $M$ of $Sl_{d_0}({\mathbb R})$ by a discrete
subgroup, $M$ being endowed with a natural $T$-invariant
probability measure $\nu$;
\item the transformation $T$ is the Sinai billiard transformation.
\end{itemize}
In these situations, we prove that $\left({1\over n^{3\over 4}}
\sum_{i=1}^n\xi_{S_k}\right)_{n\ge 1}$ converges in distribution
to the random variable $\sqrt{\sum_{m\in\mathbb Z}{\mathbb E}[\xi_0\xi_m]}\Delta_1$,
where $\Delta_1$ has the limit distribution of
$\left({1\over n^{3\over 4}}
\sum_{i=1}^n\hat\xi_{S_k}\right)_{n\ge 1}$ obtained by
Kesten and Spitzer when
$(\hat\xi_m)_m$ is a sequence of independent identically distributed
random variables with null expectation and with variance 1.
Let us notice that, in our cases, $\sum_{m\in\mathbb Z}{\mathbb E}[\xi_0\xi_m]$
is well defined and is nonnegative since
it is the limit of the variance of
${1\over \sqrt{n}}\sum_{l=0}^{n-1}\xi_l$ as $n$ goes ti infinity.

We also get the same result of convergence in distribution
for the following sequence $(\xi_k)_{k\in{\mathbb Z}}$.
Let us consider the same examples of dynamical systems $(M,{\cal F},\nu,T) $.
Instead of taking $\xi_k=f\circ T^k$, we suppose that,
conditionally to $\omega\in M$, $(\xi_k)_{k\in{\mathbb Z}}$
is an independent sequence of random variables with values in
$\{-1;1\}$.
We suppose that, conditionally to $\omega\in M$, 
$\xi_k(\omega,\cdot)$ is equal to 1 with probability
$h\circ T^k(\omega)$, for some nonnegative
H\"older continuous function $h$
with expectation $1\over 2$.
This model is envisaged by Guillotin-Plantard
and Le Ny in \cite{GuillotinLeNy} for other questions
and with other hypotheses on $(M,{\cal F},\nu,T)$ and on $f$.

Moreover we generalize this to the case when 
$\xi_k$ takes $p$ values
(conditionally to $\omega\in M$, $(\xi_k)_{k\in{\mathbb Z}}$
is an independent sequence of random variables, $\xi_k$
being equal to $\theta_j$  with probability
$f_j\circ T^k(\omega)$, with $f_1+...+f_p=1$ and with $f_1,...,f_{p}$
are nonnegative H\"older continuous functions).

In section \ref{sec:sec1}, we state a general result
under technical hypotheses of decorrelation (our theorem \ref{sec:thm}).
Section \ref{sec:proof} is devoted to the proof of this result
(the idea of the proof is inspired by one step of an inductive
method of Jan \cite{janthese,janCRAS} used in \cite{FPCMP}).

In section \ref{sec:examples}, we give some applications
of our abstract theorem \ref{sec:thm}. 
We apply our theorem \ref{sec:thm} to the examples mentionned previously
(ergodic algebraic
automorphisms of the torus, diagonal transformation of
a compact quotient of $Sl_{d_0}(\mathbb R)$, billiard transformation).
The proofs of the results of section \ref{sec:examples} are done
in sections \ref{sec:examples} and \ref{sec:examplesproof}.
\section{A technical result}\label{sec:sec1}
\begin{thm}\label{sec:thm}
Let $(S_n)_{n\ge 1}$ and $(\xi_k)_{k\in{\mathbb Z}}$ be 
two sequences of random variables defined on the same probability space
$(\Omega,{\cal T},{\mathbb P})$ such that~:
\begin{enumerate}
\item $(S_n)_{n\ge 0}$ and $(\xi_k)_{k\in{\mathbb Z}}$ are 
independent one of the other;
\item $(S_n)_{n\ge 0}$ is a simple symmetric random walk on $\mathbb Z$;
\item $(\xi_k)_{k\in{\mathbb Z}}$ 
is a stationary sequence of centered random variables
admitting moments of the fourth order;
\item we have~:
$$\sum_{p\ge 0}\sqrt{1+p}\left\vert{\mathbb E}[\xi_0\xi_p]
    \right\vert<+\infty$$
$$ \mbox{and} \ \ \sup_{N\ge 1} N^{-2}\sum_{k_1,k_2,k_3,k_4=0,...,N-1}
  \left\vert{\mathbb E}[\xi_{k_1}\xi_{k_2}\xi_{k_3}\xi_{k_4}]\right\vert<+\infty.$$
\item There exists some $C>0$, some $(\varphi_{p,s})_{p,s\in{\mathbb N}}$
and some integer $r\ge 1$ such that~:
$$\forall (p,s)\in{\mathbb N}^2,\ \ \varphi_{p+1,s}\le\varphi_{p,s} \ \ \mbox{and}
\ \ \lim_{s\rightarrow +\infty}\sqrt{s}\varphi_{rs,s}=0 $$
and such that, for all integers $n_1,n_2,n_3,n_4$ with
$0\le n_1\le n_2\le n_3\le n_4$, for all real numbers
$\alpha_{n_1},...,\alpha_{n_2}$ and $\beta_{n_3},...,\beta_{n_4}$, we have~:
$$\left\vert Cov\left(e^{i\sum_{k=n_1}^{n_2}\alpha_k\xi_k},
    e^{i\sum_{k=n_3}^{n_4}\beta_k\xi_k}\right)\right\vert
    \le C\left(1+\sum_{k=n_1}^{n_2}\vert \alpha_k\vert+
      \sum_{k=n_3}^{n_4}\vert \beta_k\vert\right) \varphi_{n_3-n_2,n_4-n_3}.$$
\end{enumerate}
Then, the sequence of random variables $\left({1\over n^{3\over 4}}
\sum_{i=1}^n\xi_{S_k}\right)_{n\ge 1}$ converges in distribution to
$\sqrt{\sum_{p\in{\mathbb Z}}{\mathbb E}[\xi_0\xi_p]}\Delta_1$,
where $\Delta_1:=\int_{\mathbb R}L_1(x)\, dB_x$,
where $(B_x)_{x\in{\mathbb R}}$ and $(b_t)_{t\ge 0}$ are two 
independent standard brownian motions and $(L_t(x))_{t\ge 0}$ 
is the local time at $x$
of $(b_t)_{t\ge 0}$, i.e. $L_t(x)=\lim_{\varepsilon\downarrow 0 }
{1\over 2\varepsilon}\int_0^t{\bf 1}_{(x-\varepsilon,x+\varepsilon)}
(b_s)\, ds$.
\end{thm}
Let us notice that the point 5 of our theorem \ref{sec:thm}
is true if $(\xi_k)_{k\in\mathbb Z}$ is a stationary
sequence of random variables satisfying the following
$\alpha$-mixing condition (cf. for example \cite{Ibragimov}, lemma 1.2)~:
$$ \lim_{n\rightarrow +\infty}\sqrt{n}\alpha_n=0,\ \
   \mbox{with}\ \ \alpha_n:=\sup_{p\ge 0;\ m\ge 0}
   \sup_{A\in \sigma(\xi_{-p},...,\xi_0)}
   \sup_{B\in \sigma(\xi_{n},...,\xi_{n+m})}
   \left\vert {\mathbb P}(A\cap B)-{\mathbb P}(A){\mathbb P}( B)\right\vert.$$

\section{Applications}\label{sec:examples}
Now let us give some examples of stationary sequences
$(\xi_k)_k$ satisfying the points 3, 4 and 5 of our theorem \ref{sec:thm}.
We say that $(M,{\cal F},\nu,T)$ is an invertible 
dynamical system if $(M,{\cal F},\nu)$ is a probability space
endowed with an invertible bi-measurable transformation $T:M\rightarrow 
M$.
%
%
%
%
\begin{hypo}\label{sec:hyp}
Let us consider an invertible dynamical system $( M,{\cal F},
 \nu, T)$ such that there exists $C_0>0$, there exist two
real sequences $(\varphi_n)_{n\ge 0}$ and $(\kappa_m)_{m\ge 0}$
and, for any function $g: M\rightarrow \mathbb C$, 
there exist $K_g^{(1)}\in[0;+\infty]$ and $K_g^{(2)}\in[0;+\infty]$
such that, for all bounded functions $g,\tilde g,h,\tilde h: M
\rightarrow{\mathbb C}$~:
\begin{enumerate}
\item for all integer $n\ge 0$, we have~:
$\vert Cov_\nu(g,h\circ  T^n)\vert\le 
     c_0\left(\Vert g\Vert_\infty\Vert h\Vert_\infty+
      \Vert h\Vert_\infty K_g^{(1)} +\Vert g\Vert_\infty K_h^{(2)}\right)
     \varphi_n$;
\item for all integer $m\ge 0$, we have~: $K_{g\circ T^{-m}}^{(1)}
     \le c_0K_g^{(1)}$;
\item for all integer $m\ge 0$, and all $k=0,...,m$, 
   we have~: $K_{h\circ T^k}^{(2)}
     \le c_0K_h^{(2)}(1+\kappa_m)$;
\item we have~: $K_{g\times\tilde g}^{(1)}\le \Vert g\Vert_\infty
   K_{\tilde g}^{(1)}+\Vert \tilde g\Vert_\infty
   K_{ g}^{(1)}$;
\item we have~: $K_{h\times\tilde h}^{(2)}\le \Vert h\Vert_\infty
   K_{\tilde h}^{(2)}+\Vert \tilde h\Vert_\infty
   K_{ h}^{(2)}$;
\item the sequence $(\varphi_n)_{n\ge 0}$
is decreasing
and there exists an integer $r\ge 1$ such that~:
$\sup_{n\ge1}{n}^6(1+\kappa_n)\varphi_{nr}<+\infty. $
\end{enumerate}
\end{hypo}
For some hyperbolic or partially hyperbolic transformations,
such properties are satisfied with $K_g^{(1)}$ some Hölder
constant of $g$ along the unstable manifolds and $K_h^{(2)}$ some Hölder
constant of $h$ along the stable-central manifolds,
with $\varphi_n=\alpha^n$ for some $\alpha\in]0;1[$ and
$\kappa_m=m^\beta$ for some $\beta\ge 0$.
Let us mention, for example, the ergodic algebraic automorphisms
of the torus as well as the diagonal transformation on compact
quotient of $Sl_{d_0}(\mathbb R)$ (cf. \cite{FPSLB}).
Moreover, in the case of the Sinai
billiard transformation, these properties come from 
\cite{ChernovDolgopyat,Chernov}.
Since the earliest work of Sinai \cite{Sin70}, these billiard
systems have been studied by many authors
(let us mention \cite{BS1,BS2,BCS1,BCS2,GO}). 
More precisely, we state~:
\begin{prop}\label{sec:PRO}
Let us consider an integer $d_0\ge 2$. 
Let $(M,{\cal F},\nu,T)$ be one of the following dynamical systems~:
\begin{itemize}
\item[{\it (i)}] $M$ is the $d_0$-dimensional torus ${\mathbb T}^{d_0}
={\mathbb R}^{d_0}/{\mathbb Z}^{d_0}$ endowed with its Borel $\sigma$-algebra 
$\cal F$
and with the normalised Haar measure $\nu$ on ${\mathbb T}^{d_0}$ and
$T$ is an algebraic automorphism of ${\mathbb T}^{d_0}$ given by
a matrix $S\in Sl_{d_0}({\mathbb Z})$ the eigenvalues of which are not
root of the unity. We endow ${\mathbb T}^{d_0}$ with the metric $d$
induced by the natural metric on ${\mathbb R}^{d_0}$.
\item[{\it (ii)}] 
$M$ is a compact quotient of $Sl_{d_0}({\mathbb R})$
by a discrete subgroup $\Gamma$ of $Sl_{d_0}({\mathbb R})$~:
$M:=\{x\Gamma;\ x\in Sl_{d_0}({\mathbb R})\}$; endowed with the normalised measure
$\nu$ induced by the Haar measure on $Sl_{d_0}({\mathbb R})$.
The transformation $T$ corresponds to the multiplication on the left
by a diagonal matrix $S=diag(T_1,...,T_{d_0})\in Sl_{d_0}({\mathbb R})$
not equal to the identity and 
such that, for all $i=1,...,d_0-1$, $T_i\ge T_{i+1}>0$.
We endow $M$ with the metric $d$ induced by
a right-translations invariant riemanian metric on $SL_{d_0}({\mathbb R})$.
\item[{\it (iii)}] 
$(M,{\cal F},\nu,T)$ is the time-discrete dynamical system
given by the discrete Sinai billiard (corresponding to the
reflection times on a scatterer). We suppose that the billiard domain
is ${\cal D}:={\mathbb T}^2\setminus\left(\bigcup_{i=1}^IO_i\right)$,
where the scatterers $O_i$ are open convex subsets of ${\mathbb T}^2$, the
closures of which are pairwise disjoint and the boundaries of which 
are $C^3$ smooth with non-null curvature.
We use the parametrisation by $(r,\varphi)$ introduced
by Sinai in \cite{Sin70} and we denote by $d$ the natural corresponding metric.
\end{itemize}
Let $\eta>0$.
We can define $g\mapsto K_g^{(1)}$
and $g\mapsto K_g^{(2)}$ such that hypothesis \ref{sec:hyp} is true and such that,
for any bounded $g:M\rightarrow\mathbb C$,
$K_g^{(1)}$ and $K_g^{(2)}$ are dominated by the H\"older constant 
$C_g^{(\eta)}$ of $g$ of order $\eta$ (eventually multiplied
by some constant).

In the case {\it (iii)}, this is still true if we replace
$C_g^{(\eta)}$ by~:
$$C_g^{(\eta,m)}:=\sup_{C\in{\cal C}_m}\sup_{x,y\in C,\ x\ne y}
{\vert g(x)-g(y)\vert\over \max(d(T^k(x),T^k(y));k=-m,..,m)^\eta},$$
for some integer $m\ge 0$,
with ${\cal C}_m=\{A\cap B;\ A\in\xi_m^u,\ B\in\xi_m^s\}$ 
with $\xi_m^u$ and $\xi_m^s$ as in \cite{Chernov} (page 7).
(We recall that, for any $k=-m,...,m$, the map $T^k$ si $C^1$
on each atom of ${\cal C}_m$).
\end{prop}
{\it Proof.\/} Let $\eta>0$.
\begin{itemize}
\item In the cases {\it (i)} and {\it (ii)}, we denote
by $\Gamma^{(s,e)}$ the set of stable-central manifolds and by $\Gamma^u$
the set of unstable manifolds.
In \cite{FPSLB}, 
each $\gamma^u\in\Gamma^u$ is endowed with some metric $d^u$
and each $\gamma^{(s,e)}\in \Gamma^{(s,e)}$ is endowed with some
metric $d^{(s,e)}$ such that there exist $\tilde c_0>0$,
$\delta_0\in]0;1[$
and $\beta>0$ such that, for
any integer $n\ge 0$, for any $\gamma^u\in\Gamma^u$
and any $\gamma^{(s,e)}\in\Gamma^{(s,e)}$, we have~:
\begin{itemize}
\item For any $y,z\in\gamma^{u}$, $d^u(y,z)\ge d(y,z)$
and for any $y',z'\in\gamma^{(s,e)}$, $d^{(s,e)}(y',z')\ge d(y',z')$.
\item For any $y,z\in \gamma^{u}$, there exists $\gamma^{u}_{(n)}$
such that $T^{-n}(y)$ and $T^{-n}(z)$ belong to $\gamma^{u}_{(n)}$
and we have~:
$d^{u}(T^{-n}(y),T^{-n}(z))\le \tilde c_0 (\delta_0)^nd^u(y,z) $.
\item For any $y,z\in \gamma^{(s,e)}$, there exists $\gamma^{(s,e)}_{(n)}$
such that $T^{n}(y)$ and $T^{n}(z)$ belong to $\gamma^{(s,e)}_{(n)}$
and we have~:
$d^{(s,e)}(T^{n}(y),T^{n}(z))\le \tilde c_0 (1+n^\beta) d^{(s,e)}(y,z) $.
\end{itemize}
Let us define~:
$$K_f^{(1)}:=\sup_{\gamma^u\in\Gamma^u}
    \sup_{y,z\in \gamma^u: y\ne z}{\vert f(y)-f(z)\vert\over 
       (d^u(y,z))^{\eta}}\ \ \mbox{and}\ \
  K_f^{(2)}:=\sup_{\gamma^{(s,e)}\in\Gamma^{(s,e)}}
    \sup_{y,z\in \gamma^{(s,e)}: y\ne z}{\vert f(y)-f(z)\vert\over 
       (d^{(s,e)}(y,z))^{\eta}} .$$
Hence, the points 2, 3, 4 and 5 of hypothesis \ref{sec:hyp} are 
satisfied with $\kappa_n=n^\beta$. 
Moreover, these two quantities are less than the H\"older
constant of order $\eta$ of $f$.

In \cite{FPSLB}, the point 1 of hypothesis \ref{sec:hyp}
is proved in the particular case {\it (ii)}.
The same proof can be used in the case {\it (i)}.
We get a sequence $(\varphi_n)_n$ decreasing exponentially fast
(cf. lemme 1.3.1 in \cite{FPSLB}).
\item Let us now consider the case {\it (iii)}.
Let us consider an integer $m\ge 0$.
Let us cconsider the set $\Gamma^s$ of homogeneous stable curves
and the set $\Gamma^u$ of homogeneous unstable curves
(see \cite{Chernov} page 7 for the definition of these curves). 
We recall that there exist 
two constants $c_1>0$ and $\delta_1\in]0;1[$
such that~:
%
\begin{itemize}
\item let $y$ and $z$ belonging to the same homogeneous unstable curve.
Then, for any integer $n\ge 0$, $T^{-n}(y)$ and $T^{-n}(z)$ belong to a
same homogeneous unstable curve and we have~:
$d(T^{-n}(y),T^{-n}(z))\le c_1{\delta_1}^n.$
Moreover, for any integer $p\ge 0$, $y$ and $z$ belong to the same atom
of $\xi_p^u$. Moreover, if $y$ and $z$ belong to the same
atom of $\xi_m^s$, then  $T^m(y)$ and $T^m(z)$
belong to a same homogeneous unstable curve.
\item let $y$ and $z$ belonging to the same homogeneous stable curve.
Then, for any integer $n\ge 0$, $T^{n}(y)$ and $T^{n}(z)$ belong to a
same homogeneous stable curve and we have~:
$d(T^{n}(y),T^{n}(z))\le c_1{\delta_1}^n.$
Moreover, for any integer $p\ge 0$, $y$ and $z$ belong to the same atom
of $\xi_p^s$. Moreover, if $y$ and $z$ belong to the same
atom of $\xi_m^u$, then  $T^{-m}(y)$ and $T^{-m}(z)$
belong to a same homogeneous stable curve.
\end{itemize}
In \cite{Chernov}, for any $y,z$, Chernov defines~:
$s_+(x,y):=\min\{n\ge 0\ :\ y\not\in\xi_n^s(x)\} $
and
$ s_-(x,y):=\min\{n\ge 0\ :\ y\not\in\xi_n^u(x)\} ,$
where $\xi_n^s(x)$ (resp. $\xi_n^u(x)$) is the atom
of $\xi_n^s$ (resp. $\xi_n^u$) containing the point $x$.

%
%
Following Chernov in \cite{Chernov} (page 15), let us introduce
the following quantities~:
$$ \tilde K_f^{(1)}:=\sup_{\gamma^u\in\Gamma^u} 
\sup_{y,z\in\gamma^u\\ y\ne z}{\vert f(y)-f(z)\vert
     \over (\delta_1)^{{\eta}s_+(y,z)}}$$
and
$$\tilde K_f^{(2)}:=\sup_{\gamma^s\in\Gamma^s} 
\sup_{y,z\in\gamma^s\\ y\ne z}{\vert f(y)-f(z)\vert
     \over (\delta_1)^{{\eta}s_-(y,z)}}.$$
In the definition of \cite{Chernov}, the suprema
are taken over all unstable and stable curves instead of
homogeneous unstable and stable curves.
However, in the proofs of theorems 4.1, 4.2 and 4.3 of
\cite{Chernov}, Chernov only uses H\"older continuity on 
homogeneous stable and unstable curves.
We observe that we have~:
$\tilde K_f^{(i)}\le 2\Vert f\Vert_\infty {\delta_1}^{-\eta m} 
+K_f^{(i)},$
with~:
$$  K_f^{(1)}:=\sup_{\gamma^u\in\Gamma^u} 
\sup_{y,z\in\gamma^u; y\ne z;s_+(y,z)\ge m+1}{\vert f(y)-f(z)\vert
     \over (\delta_1)^{{\eta}s_+(y,z)}}$$
and
$$\tilde K_f^{(2)}:=\sup_{\gamma^s\in\Gamma^s} 
\sup_{y,z\in\gamma^s; y\ne z;s_-(y,z)\ge m+1}{\vert f(y)-f(z)\vert
     \over (\delta_1)^{{\eta}s_-(y,z)}}.$$
With these definitions, we have~:
$$ K_f^{(1)}\le 
(\delta_1)^{-{\eta}(m+1)}(c_1)^{\eta}C_f^{(\eta,m)}
   \ \ \mbox{and}\ \ K_f^{(2)}\le 
(\delta_1)^{-{\eta}(m+1)}(c_1)^{\eta}C_f^{(\eta,m)}.$$
Let us prove the first inequality.
Let two points $y$ and $z$ belonging to the same 
homogeneous unstable curve such that 
$s^+(y,z)\ge m+1$. Then $y':=T^{s_+(y,z)-1}(y)$ and 
$z':=T^{s_+(y,z)-1}(z)$ belong to the same homogeneous unstable curve.
Therefore, for any $k=-m,...,m$, we have~:
\begin{eqnarray*}
d(T^k(y),T^k(z))&=&d(T^{-(s_+(y,z)-1-k)}(y'),T^{-(s_+(y,z)-1-k)}(z'))\\
&\le&c_1{\delta_1}^{s_+(y,z)-1-k}\\
&\le&c_1{\delta_1}^{s_+(y,z)-(m+1)}.
\end{eqnarray*}
Hence, since $y$ and $z$ belong to the same atom of ${\cal C}_m$, we have~:
$$\vert f(y)-f(z)\vert\le  C_f^{(\eta,m)}(c_1)^\eta{\delta_1}^{\eta s_+(y,z)}
   {\delta_1}^{-\eta(m+1)}.$$
The proof of the second inequality is analogous.

Let two points $y$ and $z$.
If $y$ and $z$ belong to the same
homogeneous unstable curve, 
then, for any integer $n\ge 0$,
we have $s_+(T^{-n}(y),T^{-n}(z))=s_+(y,z)+n$.
In the same way, if $y$ and $z$ belong to the same
homogeneous stable curve,
then for any integer $n\ge 0$,
we have $s_-(T^{n}(y),T^{n}(z))=s_-(y,z)+n$.

Hence, we get points 2, 3, 4 and 5 of hypothesis \ref{sec:hyp}
with $\kappa_n=1$.
%
%
%
%
%

Moreover, Chernov establishes the existence
of $c_3>0$ and of $\alpha_3\in]0;1[$ such that, for any integer
$n\ge 0$, for any bounded $\mathbb C$-valued functions $f$ and $g$,
we have~:
$$\vert Cov(f,g\circ T^n)\vert \le
     c_3\left(\Vert f\Vert_\infty\Vert g\Vert_\infty
       +\Vert f\Vert_\infty K_g^{(2)}+\Vert g\Vert_\infty
        K_f^{(1)}\right)(\alpha_3)^n$$
(cf. theorem 4.3 in \cite{Chernov} and the remark after theorem 4.3
in \cite{Chernov}).
This gives the points 1 and 5 of our hypothesis \ref{sec:hyp},
\end{itemize}
{\it qed.}
\begin{thm}\label{sec:appli}
Let us suppose hypothesis \ref{sec:hyp}.
Let $f:M\rightarrow{\mathbb R}$ be a bounded function.
\begin{itemize}
\item[{\it (a)}] Let us suppose that $f$ is $\nu$-centered,
that $K_f^{(1)}<+\infty$ and $K_f^{(2)}<+\infty$.
We suppose that there exists
some real number $c_1>0$ such that, for any real number $\alpha$, we have~:
$K_{\exp(i\alpha f)}^{(1)}\le c_1\vert\alpha\vert$
and $K_{\exp(i\alpha f)}^{(2)}\le c_1\vert\alpha\vert$.
Then $(\xi_k:=f\circ {T}^k)_{k\in\mathbb Z}$ 
satisfies the points 3, 4 and 5 of our theorem.
\item[{\it (b)}] Let us suppose that $f$ takes its values in $[0;1]$.
Moreover let us suppose that there exists
some $c_1>0$ such that, for any
$a,b\in\mathbb C$, we have $K_{af+b}^{(1)}\le c_1\vert a\vert$
and $K_{af+b}^{(2)}\le c_1\vert a\vert$.

Let $\left(\Omega_1:=]0;1[^{\mathbb Z},{\cal F}_1:=({\cal B}(]0;1[))
^{\times \mathbb Z},
\nu_1:=\lambda^{\otimes \mathbb Z}\right)$ where $\lambda$ is 
the Lebesgue measure on $]0;1[$.
We define $(\xi_k)_{k\in\mathbb Z}$ on the product 
$(\Omega_2:=M\times\Omega_1,{\cal F}_2:={\cal F}\otimes
{\cal F}_1,\nu_2:=\nu\otimes\nu_1)$ as follows~:
$$ \xi_k(\omega,(z_m)_{m\in{\mathbb Z}}):=
   2.{\bf 1}_{\{z_k\le f\circ T^k(\omega)\}}-1. $$
Then $(\xi_k)_{k\in \mathbb Z}$ satisfies points 3, 4 and 5 of our theorem.
\item[{\it (c)}] Let us fix an integer $p\ge 2$.  
Let us fix $p$ real numbers $\theta_1,...,\theta_{p}$
(and $\theta_0:=0$) and $p$ non-negative functions
$f_1,...,f_p:M\rightarrow[0;1]$ such that
$\int_{M}(\theta_1f_1+...+\theta_pf_p)\, d\nu=0$ 
and 
$f_1+...+f_p=1$ and such that there exists $c_2>0$
such that, for all complex numbers 
$a_1,...,a_{p-1},b$, we have 
$$\max(K_{a_1f_1+...+a_{p-1}f_{p-1}+b}^{(1)}
   ,K_{a_1f_1+...+a_{p-1}f_{p-1}+b}^{(2)})\le c_2(\vert a_1\vert+...+\vert 
    a_{p-1}\vert).$$
Let $\left(\Omega_1:=]0;1[^{\mathbb Z},{\cal F}_1:=({\cal B}(]0;1[))
^{\otimes \mathbb Z},
\nu_1:=\lambda^{\otimes \mathbb Z}\right)$ where $\lambda$ is 
the Lebesgue measure on $]0;1[$.
We define $(\xi_k)_{k\in\mathbb Z}$ on the product 
$(\Omega_2:=M\times\Omega_1,{\cal F}_2:={\cal F}\otimes
{\cal F}_1,\nu_2:=\nu\otimes\nu_1)$ as follows~:
   $$ \xi_k(\omega,(z_m)_{m\in{\mathbb Z}})=
   \sum_{l=1}^p\left(\theta_l-\theta_{l-1}\right){\bf 1}_{\{z_k\le 
   \sum_{j=1}^lf_j(T^k(\omega))\}},$$
Then $(\xi_k)_{k\in \mathbb Z}$ satisfies points 3, 4 and 5 of our theorem.

\end{itemize}
\end{thm}
Let us make some comments on the point {\it (b)}.
Conditionally to $\omega\in M$, $(\tilde\xi_k(\omega,\cdot))_{k\in{\mathbb Z}}$
is a sequence of independent random variables with values
in $\{-1;1\}$ and $\tilde\xi_k(\omega,\cdot)$ is equal to 1 with probability
$f\circ T^k(\omega)$.
This model is envisaged by Guillotin-Plantard
and Le Ny in \cite{GuillotinLeNy}.

The case {\it (c)} is a generalization of the case {\it (b)}
to the case when $\tilde\xi_k$ takes $p$ values
(conditionally to $\omega\in M$, $\tilde\xi_k(\omega,\cdot)$
is equal to $\theta_j$  with probability
$f_j\circ T^k(\omega)$).

A direct consequence of proposition \ref{sec:PRO}
and of theorem \ref{sec:appli} is~:
\begin{thm}\label{sec:thmexe}
Let $(M,{\cal F},\nu,T)$ be as in proposition
\ref{sec:PRO}. Let $\eta>0$.
Let $p\ge 2$.
Let $f,f_1,...,f_p:M\rightarrow{\mathbb R}$ be $(p+1)$ 
bounded Hölder continuous function of order
$\eta$ (or, in the case {\it (iii)} of proposition \ref{sec:PRO},
we suppose that these functions are bounded and 
such that $C_f^{(\eta,m)}<+\infty$
and $\sup_{i=1,...,p}C_{f_i}^{(\eta,m)}<+\infty$
for some integer $m\ge 0$).

We suppose that $f_1,...,f_p$ are non-negative functions satisfying
$f_1+...+f_p=1$.
\begin{itemize}
\item[{\it (a)}] Let us suppose that $f$ is $\nu$-centered.
Then $(\xi_k:=f\circ {T}^k)_{k\in \mathbb Z}$ 
satisfies points 3, 4 and 5 of our theorem.
\item[{\it (b)}] Let us suppose that $f$ takes its values in $[0;1]$ and that
we have $\int_{M}f\, d\nu={1\over 2}$.

Let $\left(\Omega_1:=]0;1[^{\mathbb Z},{\cal F}_1:=
({\cal B}(]0;1[))^{\otimes \mathbb Z},
\nu_1:=\lambda^{\otimes \mathbb Z}\right)$ where $\lambda$ is 
the Lebesgue measure on $]0;1[$.
We define $(\xi_k)_{k\in\mathbb Z}$ on the product 
$(\Omega_2:=M\times\Omega_1,{\cal F}_2:={\cal F}\otimes
{\cal F}_1,\nu_2:=\nu\otimes\nu_1)$ as follows~:
$$\tilde \xi_k(\omega,(z_m)_{m\in{\mathbb Z}}):=
   2.{\bf 1}_{\{z_k\le f\circ T^k(\omega)\}}-1. $$
Then $(\xi_k)_{k\in \mathbb Z}$ satisfies points 3, 4 and 5 of our theorem.
\item[{\it (c)}] 
Let us fix $p$ real numbers $\theta_1,...,\theta_{p}$
(and $\theta_0=0$) such that
$\int_{M}(\theta_1f_1+...+\theta_pf_p)\, d\nu=0$ 
and 
Let $\left(\Omega_1:=]0;1[^{\mathbb Z},{\cal F}_1:=({\cal B}(]0;1[))
^{\times \mathbb Z},
\nu_1:=\lambda^{\otimes \mathbb Z}\right)$ where $\lambda$ is 
the Lebesgue measure on $]0;1[$.
We define $(\xi_k)_{k\in\mathbb Z}$ on the product 
$(\Omega_2:=M\times\Omega_1,{\cal F}_2:={\cal F}\otimes
{\cal F}_1,\nu_2:=\nu\otimes\nu_1)$ as follows~:
   $$ \xi_k(\omega,(z_m)_{m\in{\mathbb Z}})=
   \sum_{l=1}^p\left(\theta_l-\theta_{l-1}\right){\bf 1}_
   {\{z_k\le \sum_{j=1}^lf_j(T^k(\omega))\}},$$
Then $(\xi_k)_{k\in \mathbb Z}$ satisfies points 3, 4 and 5 of our theorem.
\end{itemize}
\end{thm}
Let us observe that, in the case {\it (iii)} of proposition \ref{sec:PRO},
we can take the function $f$ constant on each atom of ${\cal C}_m$
for some integer $m\ge 0$. For example $f={\bf 1}_{\bigcup
_{k\ge k_0} {\mathbb H}_k}-{\bf 1}_{\bigcup
_{k\ge k_0} {\mathbb H}_{-k}}$ satisfies the case {\it (a)} of theorem
\ref{sec:thmexe} for the Sinai billiard
(with the notations $k_0$ and ${\mathbb H}_k$
of \cite{Chernov} page 5). 
In the case {\it (c)} of theorem \ref{sec:thmexe},
we can take $p=3$, $\theta_1=1$, $\theta_2=-1$,
$\theta_3$=0, $f_1={\bf 1}_{\bigcup
_{k\ge k_0} {\mathbb H}_k}$, $f_2={\bf 1}_{\bigcup
_{k\ge k_0} {\mathbb H}_{-k}}$, $f_3={\bf 1}-f_1-f_2$
in the case of the Sinai billiard (with again the notations
of \cite{Chernov} page 5).
\section{Proof of theorem \ref{sec:appli}}\label{sec:examplesproof}
In the cases ${(a)}$, ${(b)}$ and ${(c)}$, it is easy to see that
$(\xi_k)_k$ is a stationary sequence of bounded random variables
\subsection{Proof of {\it (a)}}
We have~:
\begin{eqnarray*}
\sum_{p\ge 0}\sqrt{1+p}\vert{\mathbb E}[\xi_0\xi_p]\vert
&=&\sum_{p\ge 0}\sqrt{1+p}\vert{\mathbb E}_{\nu}[f.f\circ T^p]\vert\\
&\le& c_0\Vert f\Vert_\infty\left(\Vert f\Vert_\infty+K_f^{(1)}+K_f^{(2)}
\right)
       \sum_{p\ge 0}\sqrt{1+p}\varphi_{p}<+\infty.
\end{eqnarray*}
Let us consider an integer $N\ge 1$. We have~:
$${1\over N^2} \sum_{k_1,k_2,k_3,k_4=0,...,N-1}
  \left\vert{\mathbb E}[\xi_{k_1}\xi_{k_2}\xi_{k_3}\xi_{k_4}]\right\vert
  \le {24\over N^2} \sum_{0\le k_1\le k_2\le k_3\le k_4\le N-1}
  \left\vert{\mathbb E}[\xi_{k_1}\xi_{k_2}\xi_{k_3}\xi_{k_4}]\right\vert.$$
Let us consider the set $E_N^{(1)}$ of $(k_1,k_2,k_3,k_4)$
such that $0\le k_1\le k_2\le k_3\le k_4\le N-1$ and $k_4-k_3\ge N^{1\over 3}$
We have~:
\begin{eqnarray*}
 \sum_{(k_1,k_2,k_3,k_4)\in E_N^{(1)}}
  \left\vert{\mathbb E}[\xi_{k_1}\xi_{k_2}\xi_{k_3}\xi_{k_4}]\right\vert
  &=&\sum_{(k_1,k_2,k_3,k_4)\in E_N^{(1)}}
  \left\vert Cov_{\nu}\left(f\circ T^{k_1-k_3}f\circ T^{k_2-k_3}
      f,f\circ T^{k_4-k_3}\right)\right\vert\\
  &\le&c_0N^4 \left(\Vert f\Vert_\infty^4 +\Vert f\Vert_\infty^3
      (K_f^{(2)}+3c_0K_f^{(1)})
      \right)\varphi_{\lceil N^{1\over 3}\rceil}\\
  &\le&c_0N^2 \left(\Vert f\Vert_\infty^4 +\Vert f\Vert_\infty^3
      (K_f^{(2)}+3c_0K_f^{(1)})
      \right)\sup_{n\ge 1}n^6\varphi_n.
\end{eqnarray*}
Let us consider the set $E_N^{(2)}$ of $(k_1,k_2,k_3,k_4)$
such that $0\le k_1\le k_2\le k_3\le k_4\le N-1$ and $k_4-k_3< N^{1\over 3}$
and $k_3-k_2\ge r N^{1\over 3}$.
We have~:
\begin{eqnarray*}
 \sum_{(k_1,k_2,k_3,k_4)\in E_N^{(2)}}
  \left\vert Cov\left(\xi_{k_1}\xi_{k_2},\xi_{k_3}\xi_{k_4}\right)\right\vert
  &=&\sum_{(k_1,k_2,k_3,k_4)\in E_N^{(2)}}
  \left\vert Cov_{\nu}\left(f\circ T^{k_1-k_2}f,
      (f.f\circ T^{k_4-k_3})\circ T^{k_3-k_2}\right)\right\vert\\
  &\le&c_0N^4 \left(\Vert f\Vert_\infty^4 +2c_0\Vert f\Vert_\infty^3
      (K_f^{(2)}+K_f^{(1)})
      \right)(1+\kappa_{{1\over r}\left\lceil rN^{1\over 3}\right\rceil})
      \varphi_{
       \left\lceil rN^{1\over 3}\right\rceil}\\
  &\le&c_0N^2 \left(\Vert f\Vert_\infty^4 +2c_0\Vert f\Vert_\infty^3
      (K_f^{(2)}+K_f^{(1)})
      \right)\sup_{n\ge 1}n^6(1+\kappa_n)\varphi_{rn}.
\end{eqnarray*}
Moreover, we have~:
\begin{eqnarray*}
 \sum_{(k_1,k_2,k_3,k_4)\in E_N^{(2)}}
  \left\vert {\mathbb E}[\xi_{k_1}\xi_{k_2}]{\mathbb E}
  [\xi_{k_3}\xi_{k_4}]\right\vert
  &\le&
  \left(\sum_{0\le k_1\le k_2\le N-1}
  \left\vert{\mathbb E}[\xi_{k_1}\xi_{k_2}]\right\vert\right)^2\\
  &\le&\left(N\sum_{k\ge 0}
  \left\vert{\mathbb E}_{\nu}[f.f\circ T^k]\right\vert\right)^2\\
   &\le&N^2 \left(c_0\left(\Vert f\Vert_\infty^2 +\Vert f\Vert_\infty
      (K_f^{(1)}+K_f^{(2)})\right)\sum_{k\ge 0}\varphi_k\right)^2.
\end{eqnarray*}
Let us consider the set $E_N^{(3)}$ of $(k_1,k_2,k_3,k_4)$
such that $0\le k_1\le k_2\le k_3\le k_4\le N-1$ and $k_4-k_3< N^{1\over 3}$
and $k_3-k_2< r N^{1\over 3}$ and $k_2-k_1\ge r(1+r)N^{1\over 3}$.
By the same method, we get~:
\begin{eqnarray*}
 \sum_{(k_1,k_2,k_3,k_4)\in E_N^{(3)}}
  \left\vert {\mathbb E}\left[\xi_{k_1}\xi_{k_2}\xi_{k_3}\xi_{k_4}\right]
  \right\vert
  &\le&N^2{c_0\over (1+r)^6}
     \left(\Vert f\Vert_\infty^4 +3c_0\Vert f\Vert_\infty^3
      (K_f^{(2)}+K_f^{(1)})
      \right)\sup_{n\ge 1}n^6(1+\kappa_n)\varphi_{rn}.
\end{eqnarray*}
Since the number of $(k_1,k_2,k_3,k_4)$ such that 
$0\le k_1\le k_2\le k_3\le k_4\le N-1$ and that do not belong
to $E_N^{(1)}\cup E_N^{(2)}\cup E_N^{(3)}$ is bounded by
$N^22(r+1)^3$, we get~:
$$\sup_{N\ge 1}{1\over N^2} \sum_{k_1,k_2,k_3,k_4=0,...,N-1}
  \left\vert{\mathbb E}[\xi_{k_1}\xi_{k_2}\xi_{k_3}\xi_{k_4}]\right\vert
   <+\infty.$$
Now, let us prove the point 5.
Let $n_1$, $n_2$, $n_3$ and $n_4$ be four integers such that
$0\le n_1\le n_2\le n_3\le n_4$. Let us consider any real numbers
$\alpha_{n_1},..,\alpha_{n_2}$ and $\beta_{n_3},...,\beta_{n_4}$.
We have~:

$\displaystyle 
\left\vert Cov\left(e^{i\sum_{k=n_1}^{n_2}\alpha_k\xi_k},
    e^{i\sum_{k=n_3}^{n_4}\beta_k\xi_k}\right)\right\vert
    =\left\vert Cov_{\nu}\left(e^{i\sum_{k=n_1}^{n_2}\alpha_k f\circ T^{-(n_2-k)}},
    \left(e^{i\sum_{k=n_3}^{n_4}\beta_k f\circ T^{k-n_3}}\right)
     \circ T^{n_3-n_2}\right)
    \right\vert$
\begin{eqnarray*}
    &\le&c_0\left(1+
       K_{\exp\left(i\sum_{k=n_1}^{n_2}\alpha_k f\circ T^{-(n_2-k)}\right)}
       ^{(1)} + K_{\exp\left(i\sum_{k=n_3}^{n_4}\beta_k f\circ T^{k-n_3}
       \right)}^{(2)}\right)
     \varphi_{n_3-n_2}\\
   &\le&c_0\left(1+\sum_{k=n_1}^{n_2}K_{\exp(i \alpha_k f\circ T^{-(n_2-k)}
      )}^{(1)} + \sum_{k=n_3}^{n_4}K_{\exp(i \beta_k f\circ T^{k-n_3})
      }^{(2)}\right)
     \varphi_{n_3-n_2}\\
   &\le&c_0\left(1+\sum_{k=n_1}^{n_2}c_0 c_1\vert \alpha_k \vert
      + \sum_{k=n_3}^{n_4}c_0c_1\vert \beta_k\vert (1+\kappa_{n_4-n_3}) 
      \right)
     \varphi_{n_3-n_2}.
\end{eqnarray*}
We conclude by taking $\varphi_{p,s}:=(1+\kappa_s)\varphi_p$.

\subsection{Proof of {\it (b)} and of {\it (c)}}
Let us consider {\it (c)} which is an extension
of the case {\it (b)} (by taking $p=2$, $\theta_1=1$, 
$\theta_2=-1$, $f_1=f$ and $f_2=1-f$).
Let us define the function $g:=\sum_{j=1}^p\theta_jf_j$
(in the case {\it (b)}, we have~: $g=2f-1$). This function is $\nu$-centered
and satisfies $K_g^{(1)}+K_g^{(2)}<+\infty$.
We observe that, conditionally to $\omega\in M$, the expectation
of $\xi_k(\omega,\cdot)$ is equal to $g\circ T^k(\omega)$.
Using the Fubini theorem and starting by integrating over $\Omega_1$, we 
observe that, for any integers $k$ and $l$, we have~:
${\mathbb E}[\xi_k\xi_l]={\mathbb E}_{\nu}[g\circ T^k.g\circ T^l] $
and that, for any integers $k_1,k_2,k_3,k_4$, we have~:
${\mathbb E}\left[\xi_{k_1}\xi_{k_2}\xi_{k_3}\xi_{k_4}\right]=
{\mathbb E}_{\nu}\left[\prod_{j=1}^4 g\circ T^{k_j}\right] $.
Hence, we can prove the point 4 of theorem \ref{sec:thm} as we proved it
for {\it (a)}.

Now, let us prove the point 5 of theorem \ref{sec:thm}.
We observe that, conditionally to $\omega\in M$,
the expectation of $\exp(iu\xi_k(\omega,\cdot))$ 
is $h_u\circ T^k$ with
$(h_u:=\sum_{l=1}^{p}e^{i\theta_lu}f_l)$.
This function can be rewritten~: 
$h_u=e^{i\theta_pu}+
\sum_{l=1}^{p-1}\left(e^{i\theta_lu}-e^{i\theta_pu}\right)f_l$.
The modulus of this function is bounded by 1
and we have~:
$$\max\left(K_{h_u}^{(1)},K_{h_u}^{(2)}\right)
     \le c_2 2p\max_{j=0,...,p}\vert\theta_j\vert\ \vert u\vert.$$
Let $n_1$, $n_2$, $n_3$ and $n_4$ be four integers such that
$0\le n_1\le n_2\le n_3\le n_4$. Let us consider any real numbers
$\alpha_{n_1},..,\alpha_{n_2}$ and $\beta_{n_3},...,\beta_{n_4}$.
We have~:

$\displaystyle 
\left\vert Cov\left(e^{i\sum_{k=n_1}^{n_2}\alpha_k\xi_k},
    e^{i\sum_{k=n_3}^{n_4}\beta_k\xi_k}\right)\right\vert=$
\begin{eqnarray*}
    &=&\left\vert Cov_{\nu}\left(\prod_{k=n_1}^{n_2}
    h_{\alpha_k}\circ T^k,
    \prod_{k=n_3}^{n_4}h_{\beta_k}\circ T^k
    \right)\right\vert\\
  &\le& c_0\left(1
     +c_0c_2 2p\max_{j=0,...,p}\vert\theta_j\vert\
     \left(\sum_{k=n_1}^{n_2}\vert\alpha_k\vert+
     \sum_{k=n_3}^{n_4}\vert\beta_k\vert\right)\right)
     (1+\kappa_{n_4-n_3}) 
     \varphi_{n_3-n_2}.
\end{eqnarray*}

\section{Proof of theorem \ref{sec:thm}}\label{sec:proof}
To prove our result of convergence in distribution, we use
characteristic functions.
Let us fix some real number $t$.
We will show that~:
$$\lim_{n\rightarrow +\infty}{\mathbb E}\left[\exp\left(
  {it\over n^{3\over 4}}\sum_{k=1}^n\xi_{S_k}\right)\right]
  = {\mathbb E}\left[\exp\left(it\sqrt{\sum_{p\in{\mathbb Z}}{\mathbb E}[\xi_0\xi_p]}
  \Delta_1\right)\right].$$
Let us notice that we have (cf \cite{KestenSpitzer} lemma 5, for example)~:
$${\mathbb E}\left[\exp\left(iu
  \Delta_1\right)\right]=
  {\mathbb E}\left[\exp\left(-{u^2\over 2}\int_{\mathbb R}(L_1(x))^2\, dx\right)
    \right].$$
Hence, it is enough to prove that~:
$$\lim_{n\rightarrow +\infty}{\mathbb E}\left[\exp\left(
  {it\over n^{3\over 4}}\sum_{k=1}^n\xi_{S_k}\right)\right]
  = {\mathbb E}\left[\exp\left(-{t^2\over 2}
      \sum_{p\in{\mathbb Z}}{\mathbb E}[\xi_0\xi_p]
    \int_{\mathbb R}(L_1(x))^2\, dx\right)\right].$$
In the following, for any integer $m\ge 1$ and any integer $k$, we define~:
$$ N_m(k):=Card\{j=1,...,m\ :\ S_j=k\}.$$
We notice that, for any integer $n\ge 1$, we have~:
$$\sum_{j=1}^n\xi_{S_j}=\sum_{k\in{\mathbb Z}}\xi_kN_n(k).$$
In the step 1 of our proof, we will use the following facts~:
$$C_0:=\sup_{n\ge 1}\sup_{K>0}K^2n^{-1}{\mathbb P}\left(
       \max_{m=1,...,n}\vert S_m\vert\ge K\right) <+\infty,$$
$$C_1:=\sup_{n\ge 1}\sup_{k\in\mathbb Z}n^{-{1\over 2}}\Vert N_n(k)\Vert_6<+\infty ,$$
$$C_2:=\sup_{n\ge 1}\sup_{k,\ell\in\mathbb Z}
 {\Vert N_n(\ell)-N_n(k)\Vert_2\over\sqrt{1+\vert\ell-k\vert}n^{1\over 4}}
 <+\infty.$$
The first fact comes from the Kolmogorov inequality.
We refer to \cite{KestenSpitzer} lemmas 1, 2, 3 and 4 for the proof of the other
facts.
\subsection{Step 1~: Technical part}
This is the big part of our proof.
In this part, we prove that the following quantity goes to zero
as $n$ goes to $+\infty$~:
$$\left\vert{\mathbb E}\left[\exp\left(
  {it\over n^{3\over 4}}\sum_{\ell\in\mathbb Z}\xi_{\ell}N_n(\ell)\right)\right]
     -{\mathbb E}\left[\exp\left(
  -{t^2\over 2n^{3\over 2}}\sum_{\ell,k\in{\mathbb Z}}{\mathbb E}[\xi_{\ell}
     \xi_k]N_n(\ell)^2\right)\right]\right\vert.$$
Let us fix $\varepsilon>0$. We will prove that, if $n$ is large enough,
this quantity is less than $\varepsilon$.

Our proof is inspired by a method used by Jan to establish
central limit theorem with rate of convergence
(cf. \cite{janCRAS}, \cite{janthese}, method also used in
\cite{FPCMP}).
More precisely, we adapt the idea of 
the first step of the inductive method of Jan.
\begin{itemize}
\item 
For any $K\ge 1$ and any integer $n\ge 1$, we have~:
$${\mathbb P}\left(\max_{m=1,...,n}\vert S_m\vert
    \ge K\sqrt{n}\right)\le {C_0n\over K^2n}={C_0\over K^2}. $$
Let us fix $K\ge 1$ such that $2{C_0\over K^2}< {\varepsilon\over 10}$.
Then, we have
\begin{equation}
\left\vert{\mathbb E}\left[\exp\left(
  {it\over n^{3\over 4}}\sum_{\ell\in\mathbb Z}\xi_{\ell}N_n(\ell)\right)\right]-
  {\mathbb E}\left[\exp\left(
  {it\over n^{3\over 4}}\sum_{\ell=-\lceil K\sqrt{n}\rceil}
  ^{\lceil K\sqrt{n}\rceil}\xi_{\ell}N_n(\ell)\right)\right]
  \right\vert\le 2{C_0\over K^2}< {\varepsilon\over 10}
\end{equation}
and~:
\begin{equation}
\left\vert{\mathbb E}\left[\exp\left(
  -{t^2\over 2n^{3\over 2}}\sum_{\ell,k\in{\mathbb Z}}{\mathbb E}[\xi_{\ell}
     \xi_k]N_n(\ell)^2\right)\right]-
     {\mathbb E}\left[\exp\left(
  -{t^2\over 2n^{3\over 2}}\sum_{\ell=-\lceil K\sqrt{n}\rceil}
      ^{\lceil K\sqrt{n}\rceil}\sum_{k\in\mathbb Z}{\mathbb E}[\xi_{\ell}
     \xi_k]N_n(\ell)^2\right)\right]
     \right\vert< {\varepsilon\over 10}.
\end{equation}
Hence we have to estimate~:
\begin{equation}
A_n:=\left\vert{\mathbb E}\left[\exp\left(
  {it\over n^{3\over 4}}\sum_{\ell=-\lceil K\sqrt{n}\rceil}
  ^{\lceil K\sqrt{n}\rceil}\xi_{\ell}N_n(\ell)\right)\right]
     -{\mathbb E}\left[\exp\left(
  -{t^2\over 2n^{3\over 2}}\sum_{\ell=-\lceil K\sqrt{n}\rceil}
    ^{\lceil K\sqrt{n}\rceil}\sum_{k\in\mathbb Z}{\mathbb E}[\xi_{\ell}
     \xi_k]N_n(\ell)^2\right)\right]\right\vert.
\end{equation}
\item 
In the following, $L$ will be some real number bigger than 8 and large 
enough and $n$ any integer bigger than 1 and large 
enough such that~: $ {2 K\sqrt{n}\over L}\ge L$.
We will have~:
${K\sqrt{n}\over L}\le \left\lfloor{2\lceil K\sqrt{n}\rceil+1\over L}
\right\rfloor\le {5K\sqrt{n}\over L}$.
\item 
We split our sums $\sum_{\ell=-\lceil K\sqrt{n}\rceil}
    ^{\lceil K\sqrt{n}\rceil}$
in $L$ sums over $\left\lfloor{2\lceil K\sqrt{n}\rceil+1\over L}
\right\rfloor$ terms and one sum over less than $L$ terms and so
over less than $\left\lfloor{2\lceil K\sqrt{n}\rceil+1\over L}
\right\rfloor$ terms.

For any $k=0,...,L-1$, we define~:
$$a_{k,n,L}=\exp\left(-{t^2\over 2n^{3\over 2}}
\sum_{\ell=-\lceil K\sqrt{n}\rceil+k\left\lfloor{2\lceil K\sqrt{n}
\rceil+1\over L}
\right\rfloor}^{-\lceil K\sqrt{n}\rceil+(k+1)\left\lfloor{2\lceil K\sqrt{n}
\rceil+1\over L}
\right\rfloor-1}\sum_{k\in\mathbb Z}{\mathbb E}[\xi_{\ell}
     \xi_k]N_n(\ell)^2\right)$$
and
$$ b_{k,n,L}=
\exp\left({it\over n^{3\over 4}}\sum_{\ell=-\lceil K\sqrt{n}\rceil+k\left\lfloor{2\lceil K\sqrt{n}
\rceil+1\over L}
\right\rfloor}^{-\lceil K\sqrt{n}\rceil+(k+1)\left\lfloor{2\lceil K\sqrt{n}
\rceil+1\over L}
\right\rfloor-1}\xi_\ell N_n(\ell)\right).$$
Moreover, we define~:
$$a_{L,n,L}=\exp\left(-{t^2\over 2n^{3\over 2}}
\sum_{\ell=-\lceil K\sqrt{n}\rceil+L\left\lfloor{2\lceil K\sqrt{n}
\rceil+1\over L}
\right\rfloor}^{\lceil K\sqrt{n}\rceil}\sum_{k\in\mathbb Z}{\mathbb E}[\xi_{\ell}
     \xi_k]N_n(\ell)^2\right)$$
and
$$ b_{L,n,L}=
\exp\left({it\over n^{3\over 4}}\sum_{\ell=-\lceil K\sqrt{n}\rceil+L\left\lfloor{2\lceil K\sqrt{n}
\rceil+1\over L}
\right\rfloor}^{\lceil K\sqrt{n}\rceil}\xi_\ell N_n(\ell)\right).$$
Let us notice that, for any $k=0,...,L$, we have~:
$$\vert a_{k,n,L}\vert\le 1 \ \ \mbox{and}\ \  \vert b_{k,n,L}\vert\le 1.$$
We have~:
\begin{eqnarray}
\vert A_n\vert&=& \left\vert{\mathbb E}\left[
   \prod_{k=0}^Lb_{k,n,L}-\prod_{k=0}^La_{k,n,L}\right]\right\vert \nonumber\\
   &=&\left\vert\sum_{k=0}^L{\mathbb E}\left[\left(
      \prod_{m=0}^{k-1}b_{m,n,L}\right)(b_{k,n,L}-a_{k,n,L})
      \prod_{m'=k+1}^L a_{m',n,L}\right]
      \right\vert.
\end{eqnarray}
\item
Now we explain how we can
restrict our study to the sum over the $k$ such that $(r+1)^3\le k
\le L-1$.
Indeed, the number of $k$ that do not satisfy this is equal to
$(r+1)^3+1$. Let us consider any $k=0,...,L$. We have~:
$${\mathbb E}\left[\vert b_{k,n,L}-1\vert\right] 
    \le {\vert t\vert\over n^{3\over 4}}
     {\mathbb E}\left[\left\vert\sum_{\ell=\cdots}^{\cdots}\xi_\ell
     N_n(\ell)\right\vert\right]$$
and~:
$${\mathbb E}\left[\vert a_{k,n,L}-1\vert\right] 
    \le {t^2\over 2 n^{3\over 2}}
     {\mathbb E}\left[\left\vert\sum_{\ell=\cdots}^{\cdots}
     \sum_{m\in\mathbb Z}{\mathbb E}[\xi_{\ell}
     \xi_m]N_n(\ell)^2\right\vert\right].$$
But, for any integers $\alpha$ and $\beta$ with $\beta\ge 1$, 
we have~:
\begin{eqnarray*}
{\mathbb E}\left[\left(\sum_{\ell=\alpha+1}^{\alpha+\beta}
     \xi_{\ell}N_n(\ell)\right)^2\right] &\le& \sum_{\ell=\alpha+1}
     ^{\alpha+\beta}\sum_{m=\alpha+1}^{\alpha+\beta}\vert
         {\mathbb E}[\xi_\ell\xi_m]\vert\ 
          \vert{\mathbb E}[N_n(\ell)N_n(m)]\vert\\
  &\le& \sum_{\ell=\alpha+1}
     ^{\alpha+\beta}\sum_{m=\alpha+1}^{\alpha+\beta}
        \vert{\mathbb E}[\xi_\ell\xi_m]\vert\ 
          \Vert N_n(\ell)\Vert_2\Vert N_n(m)\Vert_2\\
  &\le& (C_1)^2 n \beta\sum_{m\in\mathbb Z}\vert{\mathbb E}[\xi_0\xi_m]\vert.
\end{eqnarray*}
From which, we get~:
\begin{eqnarray*}
{\mathbb E}\left[\vert b_{k,n,L}-1\vert\right] 
    &\le& {\vert t\vert\over n^{3\over 4}}
    \sqrt{(C_1)^2 n{5 K\sqrt{n}\over L}
   \sum_{m\in\mathbb Z}\vert{\mathbb E}[\xi_0\xi_m]\vert}\\
   &\le& {\vert t\vert\over \sqrt{L}} C_1\sqrt{5 K
     \sum_{m\in\mathbb Z}\vert{\mathbb E}[\xi_0\xi_m]\vert}.
\end{eqnarray*}
Moreover we have~:
\begin{eqnarray*}
{\mathbb E}\left[\vert a_{k,n,L}-1\vert\right] 
    &\le& {t^2\over 2 n^{3\over 2}}
     {\mathbb E}\left[\sum_{\ell=\cdots}^{\cdots}
     \sum_{m\in\mathbb Z}{\mathbb E}[\xi_{\ell}
     \xi_m] N_n(\ell)^2\right]\\
 &\le&{t^2\over 2 n^{3\over 2}}{5 K\sqrt{n}\over L} 
          \sum_{m\in\mathbb Z}{\mathbb E}[\xi_{0}
     \xi_m] (C_1)^2 n\\
 &\le&{5t^2\over 2}{ K(C_1)^2\over L} 
          \sum_{m\in\mathbb Z}{\mathbb E}[\xi_{0}
     \xi_m].
 \end{eqnarray*}
Let $L_1\ge 8$ be such that for all $L\ge L_1$, we have~:
$$((r+1)^3+1){\vert t\vert\over \sqrt{L}} C_1\sqrt{5K
     \sum_{m\in\mathbb Z}\vert{\mathbb E}[\xi_0\xi_m]\vert}<{\varepsilon\over 20}$$
and
$$ ((r+1)^3+1){5t^2\over 2}{ K(C_1)^2\over L} 
          \sum_{m\in\mathbb Z}{\mathbb E}[\xi_{0}
     \xi_m] <{\varepsilon\over 20}.$$
Then, if we have $L\ge L_1$ and $n\ge 1$ such that $ 
{2 K\sqrt{n}\over L}\ge L$, we have~:
\begin{equation}
{\mathbb E}\left[\vert b_{L,n,L}-a_{L,n,L}\vert\right]+
\sum_{k=0}^{(r+1)^3-1}{\mathbb E}\left[\vert b_{k,n,L}-a_{k,n,L}\vert\right]
<{\varepsilon\over 10}.
\end{equation}
It remains to estimate~:
\begin{equation}
\sum_{k=(r+1)^3}^{L-1}\left\vert{\mathbb E}\left[\left(
      \prod_{m=0}^{k-1}b_{m,n,L}\right)(b_{k,n,L}-a_{k,n,L})
      \prod_{m'=k+1}^L a_{m',n,L}\right]
      \right\vert.
\end{equation}
\item We estimate~:
$$B_{n,L}:=\sum_{k=(r+1)^3}^{L-1}\left\vert{\mathbb E}\left[\left(
      \prod_{m=0}^{k-(r+1)^3}b_{m,n,L}\right)
      \left(\left(\prod_{m=k-(r+1)^{3}+1}^{k-(r+1)^{2}}b_{m,n,L}
      \right)-1\right)\times\right.\right.$$
$$\left.\left.\times\left(\left(\prod_{m=k-(r+1)^{2}+1}^{k-r-1}b_{m,n,L}
      \right)-1\right)
\left(\prod_{m'=k-r}^{k-1}b_{m',n,L}\right)
  (b_{k,n,L}-a_{k,n,L})\prod_{m'=k+1}^L 
  a_{m',n,L}\right] \right\vert.$$
We have~:
$$B_{n,L}\le \sum_{k=(r+1)^3}^{L-1}
\left\Vert \left(\prod_{m=k-(r+1)^{3}+1}^{k-(r+1)^{2}}b_{m,n,L}
      \right)-1\right\Vert_3
      \left\Vert \left(\prod_{m=k-(r+1)^{2}+1}^{k-r-1}b_{m,n,L}
      \right)-1\right\Vert_3\left\Vert b_{k,n,L}-a_{k,n,L}\right\Vert_3.$$
\begin{itemize}
\item We have~:
$$\Vert b_{k,n,L}-1\Vert_3\le {\vert t\vert\over n^{3\over 4}}
   \left\Vert \sum_{\ell=\cdots}^{\cdots}\xi_\ell N_n(\ell)\right\Vert_3. $$
For any integers $\alpha$ and $\beta$ with $\beta\ge 1$, we have~:
\begin{eqnarray}
{\mathbb E}\left[\left(\sum_{\ell=\alpha+1}^{\alpha+\beta}\xi_\ell N_n(\ell)
\right)^4\right]
&\le&\sum_{\ell_1,\ell_2,\ell_3,\ell_4=\alpha+1}^{\alpha+\beta}
\left\vert{\mathbb E}\left[
    \xi_{\ell_1}\xi_{\ell_2}\xi_{\ell_3}\xi_{\ell_4}\right]\right\vert
      (C_1)^4 n^2\nonumber\\
&\le&(C_1)^4 n^2 C'_2\beta^2 \label{sec:norme4}.
\end{eqnarray}
with $C'_2:=\sup_{N\ge 1} N^{-2}\sum_{k_1,k_2,k_3,k_4=0,...,N-1}
  \left\vert{\mathbb E}[\xi_{k_1}\xi_{k_2}\xi_{k_3}\xi_{k_4}]\right\vert$.
Hence, we have~:
\begin{eqnarray*}
\Vert b_{k,n,L}-1\Vert_3&\le& {\vert t\vert\over n^{3\over 4}}
   \left((C_1)^4 n^2 C'_2\left( {5K\sqrt{n}\over L}\right)^2\right)^{1\over 4}\\
    &\le& {\vert t\vert}C_1
   \left(  C'_2\right)^{1\over 4}\sqrt{5K\over L}.
\end{eqnarray*}
\item We have~:
\begin{eqnarray*}
\Vert a_{k,n,L}-1\Vert_3&\le& { t^2\over 2n^{3\over 2}}
\sum_{k\in\mathbb Z}{\mathbb E}[\xi_{0}
     \xi_k] \left\Vert \sum_{\ell=\cdots}^{\cdots}
      N_n(\ell)^2\right\Vert_3\\
&\le&{ t^2\over 2n^{3\over 2}}
\sum_{k\in\mathbb Z}{\mathbb E}[\xi_{0}
     \xi_k] \sum_{\ell=\cdots}^{\cdots}\left\Vert 
      N_n(\ell)\right\Vert_6^2\\
&\le&{ t^2\over 2n^{3\over 2}}
\sum_{k\in\mathbb Z}{\mathbb E}[\xi_{0}
     \xi_k] {5K\sqrt{n}\over L}(C_1)^2 n\\
&\le&{ 5t^2\over 2}
\sum_{k\in\mathbb Z}{\mathbb E}[\xi_{0}
     \xi_k] {K\over L}(C_1)^2.
\end{eqnarray*}
\item Using formula (\ref{sec:norme4}),we get~:
\begin{eqnarray*}
\left\Vert \left(\prod_{m=k-(r+1)^{3}+1}^{k-(r+1)^{2}}b_{m,n,L}
      \right)-1\right\Vert_3
&\le& {\vert t\vert\over n^{3\over 4}}  
   \left\Vert \sum_{\ell=-\lceil K\sqrt{n}\rceil
      +(k-(r+1)^3+1)\left\lfloor{2\lceil K\sqrt{n}
\rceil+1\over L}
\right\rfloor}^{-\lceil K\sqrt{n}\rceil
      +(k-(r+1)^2+1)\left\lfloor{2\lceil K\sqrt{n}
\rceil+1\over L}
\right\rfloor-1}\xi_\ell N_n(\ell)\right\Vert_3\\
&\le& {\vert t\vert\over n^{3\over 4}}  \left( (C_1)^4n^2C'_2
   (r(r+1)^2)^2\left({5 K\sqrt{n}\over L}\right)^2\right)^{1\over 4}\\
&\le& {\vert t\vert}  C_1 (C'_2)^{1\over 4}
   \sqrt{r}(r+1)\sqrt{5 K\over {L}}.
\end{eqnarray*}
\item Analogously, we get~:
$$\left\Vert \left(\prod_{m=k-(r+1)^{2}+1}^{k-r-1}b_{m,n,L}
      \right)-1\right\Vert_3
\le {\vert t\vert} C_1 (C'_2)^{1\over 4}
   \sqrt{r}(r+1)\sqrt{5 K\over L}.$$
\end{itemize}
Hence, we have~:
\begin{eqnarray*}
B_{n,L}&\le& L   \left({\vert t\vert} C_1(C'_2)^{1\over 4}
   \sqrt{r}(r+1)\sqrt{5 K\over L}\right)^2
   \left({\vert t\vert} C_1
   \left(  C'_2\right)^{1\over 4}\sqrt{5K\over L}+
   { 5t^2\over 2}
\sum_{k\in\mathbb Z}{\mathbb E}[\xi_{0}
     \xi_k] {K\over L}(C_1)^2
   \right)\\
&\le&   {\vert t\vert}^2(C_1)^2(C'_2)^{1\over 2}
   r(r+1)^2{5 K}
   \left({\vert t\vert}C_1
   \left( C'_2\right)^{1\over 4}\sqrt{5K\over L}+
   { 5t^2\over 2}
\sum_{k\in\mathbb Z}{\mathbb E}[\xi_{0}
     \xi_k] {K\over L}(C_1)^2
   \right).
\end{eqnarray*}
Let $L'_1\ge L_1$ be such that, for all $L\ge L_1$, 
the right term of this last inequality is less than 
$\varepsilon\over 10$.

Then, for any $L\ge L'_1$ and any $n\ge 1$ such that
$ {2 K\sqrt{n}\over L}\ge L$, we have~:
$B_{n,L}\le {\varepsilon\over 10}$.
\item In the following, we suppose $L\ge L'_1$ and 
$ {2 K\sqrt{n}\over L}\ge L$.
It remains to estimate~:
$$\sum_{k=(r+1)^3+1}^{L-1} C_{n,k,L,1,3}+C_{n,k,L,1,2}+C_{n,k,L,2,3} $$
where $C_{n,k,L,j_0,j_1}$ is the following quantity~:
$$\left\vert{\mathbb E}\left[\left(
      \prod_{m=0}^{k-(r+1)^{j_1}}b_{m,n,L}\right)
      \left(\prod_{m=k-(r+1)^{j_0}+1}^{k-1}b_{m,n,L}
      \right)(b_{k,n,L}-a_{k,n,L})\prod_{m'=k+1}^L 
  a_{m',n,L}\right] \right\vert$$
\item
Let $j_0,j_1$ be fixed. We estimate $C_{n,k,L,j_0,j_1}$. We have~:
$$C_{n,k,L,j_0,j_1}\le D_{n,k,L,j_0,j_1}+ E_{n,k,L,j_0,j_1},$$
where~:
$$\displaystyle D_{n,k,L,j_0,j_1}:=
  \left\vert{\mathbb E}\left[Cov_{\vert(S_p)_p}  \left(\Delta_{n,k,L,j_1},
       \Gamma_{n,k,L,j_0}\right)\prod_{m'=k+1}^L 
  a_{m',n,L}\right]\right\vert$$

and
$$E_{n,k,L,j_0,j_1}:
  \left\vert{\mathbb E}\left[{\mathbb E}\left[\left. \Delta_{n,k,L,j_1} 
      \right\vert{(S_p)_p}\right]
      {\mathbb E}\left[\left. \Gamma_{n,k,L,j_0}
      \right\vert(S_p)_p\right]\prod_{m'=k+1}^L 
  a_{m',n,L} \right] \right\vert.$$
with $\Delta_{n,k,L,j_1}:=\prod_{m=0}^{k-(r+1)^{j_1}}b_{m,n,L}$
and $\Gamma_{n,k,L,j_0}:=\left(\prod_{m=k-(r+1)^{j_0}+1}^{k-1}b_{m,n,L}
      \right)(b_{k,n,L}-a_{k,n,L})$.
\item 
Control of the terms with the product of the expectations.

Let $j_0,j_1$ be fixed. Let $k=(r+1)^3,...,L-1$.
We can notice that $E_{n,k,L,j_0,j_1}$ is bounded from away
by the following quantity~:
$$ F_{n,k,L,j_0,j_1}:=
  {\mathbb E}\left[\left\vert
      {\mathbb E}\left[\left.
      \prod_{m=k-(r+1)^{j_0}+1}^{k}b_{m,n,L}-
      \left(\prod_{m=k-(r+1)^{j_0}+1}^{k-1}b_{m,n,L}
      \right)a_{k,n,L}
      \right\vert (S_p)_p\right]\right\vert 
  \right].$$
We use the Taylor expansions of the exponential function.
To simplify expressions, we will use the following notation~:
$$\forall m\ge 0,\ \ 
\alpha_{(m)}:= -\lceil K\sqrt{n}\rceil+m\left\lfloor{2\lceil K\sqrt{n}
\rceil+1\over L}\right\rfloor.$$
\begin{itemize}
\item 
Let us show that, in $F_{n,k,L,j_0,j_1}$, we can replace
$$\prod_{m=k-(r+1)^{j_0}+1}^{k}b_{m,n,L} =
\exp\left({it\over n^{3\over 4}}\sum_{\ell=\alpha_{(k-(r+1)^{j_0}+1)}}
  ^{\alpha_{(k+1)}-1}\xi_\ell N_n(\ell)\right)$$
by the formula given by the Taylor expansion of the exponential 
function at the
second order~:
\begin{equation}\label{sec:Taylor1}
1+{it\over n^{3\over 4}}\sum_{\ell=\alpha_{(k-(r+1)^{j_0}+1)}}
  ^{\alpha_{(k+1)}-1}\xi_\ell N_n(\ell)-
  {t^2\over 2n^{3\over 2}}\left(\sum_{\ell=\alpha_{(k-(r+1)^{j_0}+1)}}
  ^{\alpha_{(k+1)}-1}\xi_\ell N_n(\ell)\right)^2 .
\end{equation}
Indeed the $L^1$-norm of the error between these two quantities is
less than~:
$${\vert t\vert^3\over 6 n^{9\over 4}}
   {\mathbb E}\left[\left\vert \sum_{\ell=\alpha_{(k-(r+1)^{j_0}+1)}}
  ^{\alpha_{(k+1)}-1}\xi_\ell N_n(\ell)\right\vert^3\right] $$
which, according to formula (\ref{sec:norme4}), is less than~:
$${\vert t\vert^3\over 6 n^{9\over 4}}\left(
(C_1)^4n^2C'_2\left(((r+1)^{j_0}){5K\sqrt{n}\over L}\right)^2\right)
    ^{3\over 4}
    ={\vert t\vert^3\over 6 }(C_1)^3
(C'_2)^{3\over 4}\left((r+1)^{j_0}{5K\over L}\right)^{3\over 2}.
   $$
Hence, the sum over $k=(r+1)^3,...,L-1$ of the $L^1$-norm of these errors
is less than~:
$${1\over\sqrt{L}}{\vert t\vert^3\over 6 }
(C'_2)^{3\over 4}(C_1)^3\left((r+1)^{j_0}{5K}\right)^{3\over 2}.
 $$
Let us consider $L_2\ge L''_1$ such that, for all $L\ge L_2$, 
this last quantity is less than $\varepsilon\over 10$.
\item Let us introduce
$ Y_k:=\sum_{\ell=\alpha_{(k-(r+1)^{j_0}+1)}}
  ^{\alpha_{(k)}-1}\xi_\ell N_n(\ell)$ and 
$ Z_k:=\sum_{\ell=\alpha_{(k)}}^{\alpha_{(k+1)}-1}
     \sum_{m\in\mathbb Z}{\mathbb E}[\xi_\ell\xi_m]N_n(\ell)^2$.
We show that, in $F_{n,k,L,j_0,j_1}$, we can replace
$$\left(\prod_{m=k-(r+1)^{j_0}+1}^{k-1}b_{m,n,L}\right)
    a_{k,n,L} =
\exp\left({it\over n^{3\over 4}}Y_k-
  {t^2\over 2n^{3\over 2}}Z_k\right)$$
by the formula given by the Taylor expansion of the exponential 
function at the second order~:
\begin{equation}\label{sec:Taylor2}
1+{it\over n^{3\over 4}}Y_k-{t^2\over 2n^{3\over 2}}Z_k
    +{1\over 2}\left({it\over n^{3\over 4}}Y_k-{t^2\over 2n^{3\over 2}}Z_k\right)^2 ,
\end{equation}

Indeed, the $L^1$-norm of the error between these two quantities is
less than~:
$${1\over 6}{\mathbb E}\left[\left\vert 
{it\over n^{3\over 4}}Y_k-{t^2\over 2n^{3\over 2}}Z_k\right\vert^3\right]
\le {4\over 3}{\mathbb E}\left[\left\vert 
{it\over n^{3\over 4}}Y_k\right\vert^3+
\left\vert{t^2\over 2n^{3\over 2}}Z_k\right\vert^3\right].
$$
According to formula (\ref{sec:norme4}), we have~:
$$ {4\over 3}{\mathbb E}\left[\left\vert 
{it\over n^{3\over 4}}Y_k\right\vert^3\right]
  \le {4\vert t\vert^3\over 3}(C_1)^3
(C'_2)^{3\over 4}\left((r+1)^{j_0}{5K\over L}\right)^{3\over 2}.
$$
Moreover, we have~:
\begin{eqnarray*}
{4\over 3}{\mathbb E}\left[
\left\vert{t^2\over 2n^{3\over 2}}Z_k\right\vert^3\right]
&=&{t^6\over 6n^{9\over 2}}\left(\sum_{m\in\mathbb Z}{\mathbb E}
  [\xi_0\xi_m]\right)^3{\mathbb E}\left[
\sum_{\ell_1\ell_2,\ell_3=\alpha_{(k)}}^{\alpha_{(k+1)}-1}
     N_n(\ell_1)^2N_n(\ell_2)^2N_n(\ell_3)^2\right]\\
&\le&{t^6\over 6n^{9\over 2}}\left(\sum_{m\in\mathbb Z}{\mathbb E}
  [\xi_0\xi_m]\right)^3 \left({5K\sqrt{n}\over L}\right)^3
    (C_1)^6n^3.
\end{eqnarray*}
The sum over $k=(r+1)^3,...,L-1$ of the $L^1$-norm of these errors is less 
than~:
$${1\over\sqrt{L}}{4\vert t\vert^3\over 3 }
(C'_2)^{3\over 4}(C_1)^3\left((r+1)^{j_0}{5K}\right)^{3\over 2}
+{1\over L^2}{t^6\over 6}\left(\sum_{m\in\mathbb Z}{\mathbb E}
  [\xi_0\xi_m]\right)^3 \left({5K}\right)^3
    (C_1)^6. $$
Let us consider $L'_2\ge L_2$ such that, for all $L\ge L'_2$,
this last quantity is less than $\varepsilon\over 10$.
\item Now we show that, in formula (\ref{sec:Taylor2}), 
we can ommit the term with $(Z_k)^2$.
Indeed, we have~:
$${1\over 2}{\mathbb E}\left[\left({t^2\over 2n^{3\over 2}}Z_k\right)^2\right]
\le {t^4\over 8n^{3}} \left(\sum_{m\in\mathbb Z}{\mathbb E}
   [\xi_0\xi_m]\right)^2\left({5K\sqrt{n}\over L}\right)^2
   (C_1)^4n^2.$$
The sum over $k=(r+1)^3,...,L-1$ of the $L^1$-norm of these errors is less 
than~:
$${1\over L}{t^4\over 8} \left(\sum_{m\in\mathbb Z}{\mathbb E}
   [\xi_0\xi_m]\right)^2\left({5K}\right)^2
   (C_1)^4.$$
Let us consider $L''_2\ge L'_2$ such that, for all $L\ge L''_2$,
this last quantity is less than $\varepsilon\over 10$.
\item From now, we fix $L:=L''_2$ and we consider 
an integer $n\ge {L^4\over 4K^2}$.
\item Hence, it remains to estimate the following quantity called
$G_{n,k,L,j_0,j_1}$~:

$\displaystyle{\mathbb E}\left[\left\vert 
{\mathbb E}\left[  {it\over n^{3\over 4}}
(Y_k+W_k)-
  {t^2\over 2n^{3\over 2}}\left(Y_k+W_k\right)^2 
-{it\over n^{3\over 4}}Y_k+{t^2\over 2n^{3\over 2}}Z_k
+
\right.\right.\right.$
$$\left.\left.\left.\left.  +{t^2\over 2n^{3\over 2}}(Y_k)^2
    +{it\over n^{3\over 4}}Y_k{t^2\over 2n^{3\over 2}}Z_k
\right\vert (S_p)_p\right]    \right\vert\right]$$
with
$W_k:= \sum_{\ell=\alpha_{(k)}}
  ^{\alpha_{(k+1)}-1}\xi_\ell N_n(\ell)$.
Using the fact that the $\xi_k$ are centered and independent of
$(S_p)_p$, we get~:
\begin{eqnarray*}
G_{n,k,L,j_0,j_1}&=&{\mathbb E}\left[\left\vert 
{\mathbb E}\left[\left.  -
  {t^2\over 2n^{3\over 2}}\left(Y_k+W_k\right)^2 
+{t^2\over 2n^{3\over 2}}Z_k
   +{t^2\over 2n^{3\over 2}}(Y_k)^2
    \right\vert (S_p)_p\right]    \right\vert\right]\\
&=&{t^2\over 2n^{3\over 2}}{\mathbb E}\left[\left\vert 
{\mathbb E}\left[\left.  
  \left(W_k\right)^2 +2 W_kY_k-Z_k
       \right\vert (S_p)_p\right]    \right\vert\right].
\end{eqnarray*}
Let us notice that we have~:
$$Z_k:= \sum_{\ell=\alpha_{(k)}}^{\alpha_{(k+1)}-1}
     \left(
     {\mathbb E}[(\xi_\ell)^2]N_n(\ell)^2+
       2\sum_{m\le \ell -1}{\mathbb E}[\xi_\ell\xi_m]N_n(\ell)^2\right).$$
\item Let us show that, in the last expression of $G_{n,k,L,j_0,j_1}$, we
can replace $Z_k$ by~:
$$\tilde Z_k:= \sum_{\ell=\alpha_{(k)}}^{\alpha_{(k+1)}-1}
     \left(
     {\mathbb E}[(\xi_\ell)^2]N_n(\ell)^2+
       2\sum_{m\le \ell -1}{\mathbb E}[\xi_\ell\xi_m]N_n(\ell)N_n(m)\right).$$
Indeed, we have~:
\begin{eqnarray*}
{t^2\over 2n^{3\over 2}}{\mathbb E}\left[\left\vert Z_k-\tilde Z_k\right\vert\right]
&\le& {t^2\over n^{3\over 2}} \sum_{\ell=\alpha_{(k)}}^{\alpha_{(k+1)}-1}
   \sum_{m\le \ell -1}\vert{\mathbb E}[\xi_\ell\xi_m]
  \vert \ \Vert N_n(\ell)\Vert_2 \Vert N_n(m)-N_n(\ell)\Vert_2\\
&\le& {t^2\over n^{3\over 2}}{5K\sqrt{n}\over L}\sum_{p\ge  1}
  \vert{\mathbb E}[\xi_0\xi_p]
  \vert  C_1\sqrt{n}C_2n^{1\over 4}\sqrt{1+ p}.
\end{eqnarray*}     
The sum over $k=(r+1)^3,...,L-1$ of these quantities is less than~:
$${t^2\over n^{1\over 4}}{5K} C_1C_2\sum_{p\ge  1}
  \sqrt{1+ p}\vert{\mathbb E}[\xi_0\xi_p]
  \vert ,$$
which goes to zero when $n$ goes to infinity. 
Hence, there exists some $n_0\ge {L^4\over 4K^2}$ such that, for any
integer $n\ge n_0$, this sum is less than $\varepsilon\over 10$. 
\item 
Hence we have to estimate~:       
$$\tilde G_{n,k,L,j_0,j_1}={t^2\over 2n^{3\over 2}}{\mathbb E}\left[\left\vert 
{\mathbb E}\left[\left.  
  (W_k)^2+2W_kY_k
       \right\vert (S_p)_p\right]-\tilde Z_k    \right\vert\right].$$
We have~:
$$
{\mathbb E}\left[\left.  
  (W_k)^2
       \right\vert (S_p)_p\right]
       =\sum_{\ell=\alpha_{(k)}}
  ^{\alpha_{(k+1)}-1}
  \left({\mathbb E}[(\xi_{\ell})^2] (N_n(\ell))^2
  +2\sum_{m=\alpha_{(k)}}^{\ell-1}{\mathbb E}[\xi_{\ell}\xi_{m}] N_n(\ell)N_n(m)
    \right).$$
Hence we have~:
$$
{\mathbb E}\left[\left.  
  (W_k)^2+2W_kY_k
       \right\vert (S_p)_p\right]
       =\sum_{\ell=\alpha_{(k)}}
  ^{\alpha_{(k+1)}-1}
  \left({\mathbb E}[(\xi_{\ell})^2] (N_n(\ell))^2
  +2\sum_{m=\alpha_{(k-(r+1)^{j_0}+1)}}^{\ell-1}{\mathbb E}[\xi_{\ell}\xi_{m}] N_n(\ell)N_n(m)
    \right).$$
\end{itemize}
We get~:
\begin{eqnarray*}
 \tilde G_{n,k,L,j_0,j_1}
   &=& {t^2\over n^{3\over 2}}{\mathbb E}\left[\left\vert\sum_{\ell=\alpha_{(k)}}
  ^{\alpha_{(k+1)}-1}
  \sum_{m\le\alpha_{(k-(r+1)^{j_0}+1)}-1}{\mathbb E}[\xi_{\ell}\xi_{m}] 
    N_n(\ell)N_n(m)
   \right\vert\right]\\
  &\le&{t^2\over n^{3\over 2}}{5K\sqrt{n}\over L}
   \sum_{m\ge (r+1)^{j_0}{K\sqrt{n}\over L}}
         \vert{\mathbb E}[\xi_{0}\xi_{m}] \vert (C_1)^2n .
\end{eqnarray*}
The sum over $k=(r+1)^3,...,L-1$ of these quantities is less than~:
$${t^2}{5K}
   \sum_{m\ge (r+1)^{j_0}{K\sqrt{n}\over L}}
         \vert{\mathbb E}[\xi_{0}\xi_{m}] \vert (C_1)^2,$$
which goes to zero when $n$ goes to infinity.   
Hence, there exists some $n'_0\ge n_0$ such that, for any
integer $n\ge n_0$, this sum is less than $\varepsilon\over 10$. 
\item
Control of the covariance terms.

Let $j_0,j_1$ be fixed. Let $k=(r+1)^3,...,L-1$.
We have~:
\begin{eqnarray*}
 D_{n,k,L,j_0,j_1}
 &\le&\left\vert{\mathbb E}\left[Cov_{\vert(S_p)_p}  \left(
      \prod_{m=0}^{k-(r+1)^{j_1}}b_{m,n,L},
      \prod_{m=k-(r+1)^{j_0}+1}^{k}b_{m,n,L}
      \right)\prod_{m'=k+1}^L 
  a_{m',n,L}\right] \right\vert+\\
  &\ &+
  \left\vert{\mathbb E}\left[Cov_{\vert(S_p)_p}  \left(
      \prod_{m=0}^{k-(r+1)^{j_1}}b_{m,n,L},
      \prod_{m=k-(r+1)^{j_0}+1}^{k-1}b_{m,n,L}
      \right)\prod_{m'=k}^L a_{m',n,L}\right] \right\vert.
\end{eqnarray*}
But we have~:
$$\prod_{m=\alpha}^{\alpha+\beta}b_{m,n,L}
   =\exp\left({it\over n^{3\over 4}}\sum_{\ell=
     -\lceil K\sqrt{n}\rceil+\alpha
       \left\lfloor{2\lceil K\sqrt{n}\rceil +1\over L}\right\rfloor}
       ^{-\lceil K\sqrt{n}\rceil+(\alpha+\beta+1)
       \left\lfloor{2\lceil K\sqrt{n}\rceil +1\over L}\right\rfloor-1}
       \xi_\ell N_n(\ell)\right).$$
Therefore, according to point 4 of the hypothesis of our theorem, we have~:
$$D_{n,k,L,j_0,j_1}\le 2{\mathbb E}
\left[C\left(1+{\vert t\vert\over n^{3\over 4}}
\sum_{\ell\in\mathbb Z}
    N_n(\ell)\right)\left({rK\sqrt{n}\over 2L}\right)^{-{1\over 2}}
          \sup_{s\ge r{K\sqrt{n}\over 2L}}\sqrt{s}\varphi_{rs,s}
          \right].$$
Hence, we have~:
$$\sum_{k=(r+1)^3}^{L-1}D_{n,k,L,j_0,j_1}\le 2CL\sqrt{L}
C(1+\vert t\vert n^{1\over 4}){n^{-{1\over 4}}\sqrt{2}\over\sqrt{rK}}
          \sup_{s\ge r{K\sqrt{n}\over 2L}}\sqrt{s}\varphi_{rs,s}.$$

which goes to zero as $n$ goes to infinity.
Hence, there exists some $N_0\ge n'_0$ such that, for any
integer $n\ge n_0$, this sum is less than $\varepsilon\over 10$. 
\end{itemize}
Therefore, there exists $N_0$ (depending on $t$ and on $\varepsilon$) 
such that, for any integer
$n\ge N_0$, we have~:
$$\left\vert{\mathbb E}\left[\exp\left(
  {it\over n^{3\over 4}}\sum_{\ell\in\mathbb Z}\xi_{\ell}N_n(\ell)\right)\right]
     -{\mathbb E}\left[\exp\left(
  -{t^2\over 2n^{3\over 2}}\sum_{\ell,k\in{\mathbb Z}}{\mathbb E}[\xi_{\ell}
     \xi_k]N_n(\ell)^2\right)\right]\right\vert<\varepsilon.$$
This ends the step 1 of our proof.
\subsection{Step 2~: Conclusion}
In the previous section we proved that~:
$$
\lim_{n\rightarrow +\infty}\left\vert{\mathbb E}\left[\exp\left(
  {it\over n^{3\over 4}}\sum_{\ell\in\mathbb Z}\xi_{\ell}N_n(\ell)\right)\right]
     -{\mathbb E}\left[\exp\left(
  -{t^2\over 2n^{3\over 2}}\sum_{\ell,k\in{\mathbb Z}}{\mathbb E}[\xi_{\ell}
     \xi_k]N_n(\ell)^2\right)\right]\right\vert=0.$$
According to \cite{KestenSpitzer} lemma 6, we know that~:
$\left({1\over n^{3\over 2}}\sum_{\ell\in{\mathbb Z}}N_n(\ell)^2\right)
_{n\ge 1} $ converges in distribution to $Z_1:=\int_{\mathbb R}(L_1(x))^2\, dx$.
Hence, we get~:
$$\lim_{n\rightarrow +\infty}{\mathbb E}\left[\exp\left(
  {it\over n^{3\over 4}}\sum_{\ell\in\mathbb Z}\xi_{\ell}N_n(\ell)\right)\right]
={\mathbb E}\left[\exp\left(-{t^2\over 2}
    \sum_{k\in{\mathbb Z}}{\mathbb E}[\xi_{0}
     \xi_k]\int_{\mathbb R}(L_1(x))^2\, dx\right)\right]. $$

\end{document}